\documentclass{vgtc}

\ifpdf
  \pdfoutput=1\relax                   
  \usepackage{graphicx}                
  \DeclareGraphicsExtensions{.pdf,.png,.jpg,.jpeg} 
\else
  \usepackage{graphicx}                
  \DeclareGraphicsExtensions{.eps}     
\fi%

\graphicspath{{figures/}{pictures/}{images/}{./}} 

\usepackage{microtype}                 
\PassOptionsToPackage{warn}{textcomp}  
\usepackage{textcomp}                  
\usepackage{mathptmx}                  
\usepackage{times}                     
\usepackage{cite}                      
\usepackage{tabu}                      
\usepackage{booktabs}                  

\onlineid{0}

\vgtccategory{Research}

\vgtcinsertpkg

\title{Planar Symmetry Detection and Quantification using the Extended Persistent Homology Transform}

\author{Nicholas Bermingham \thanks{email: nicholas.bermingham@anu.edu.au}\\ %
        \scriptsize Australian National University %
\and Vanessa Robins\thanks{email:vanessa.robins@anu.edu.au}\\ %
     \scriptsize Australian National University %
\and Katharine Turner\thanks{email: katharine.turner@anu.edu.au}\\ %
     \parbox{1.4in}{\scriptsize \centering Australian National University}}

\teaser{
  \centering
  a)\includegraphics[scale = 0.18]{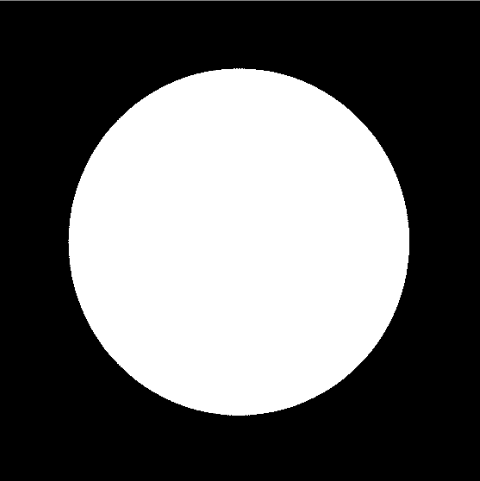}
     \includegraphics[scale = 0.2]{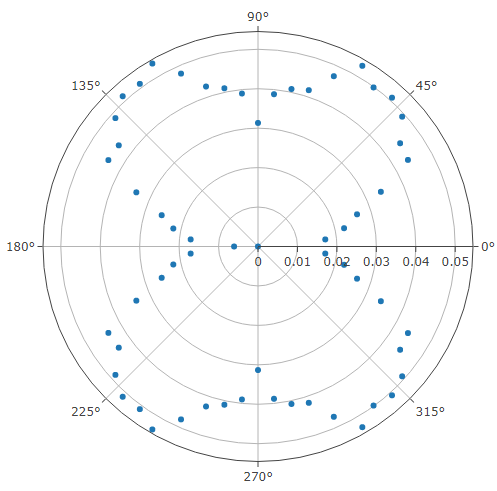}
     \includegraphics[scale = 0.2]{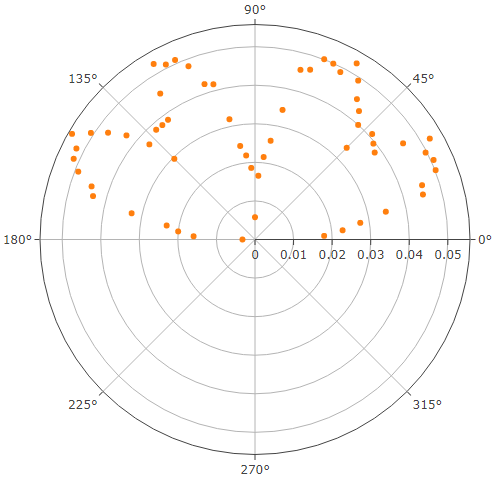}
    b)\includegraphics[scale = 0.18]{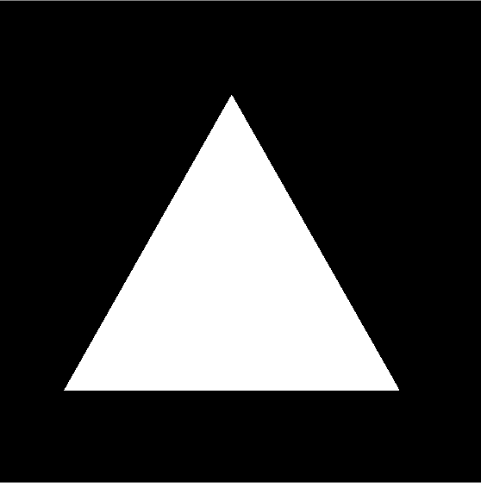}
    \includegraphics[scale = 0.2]{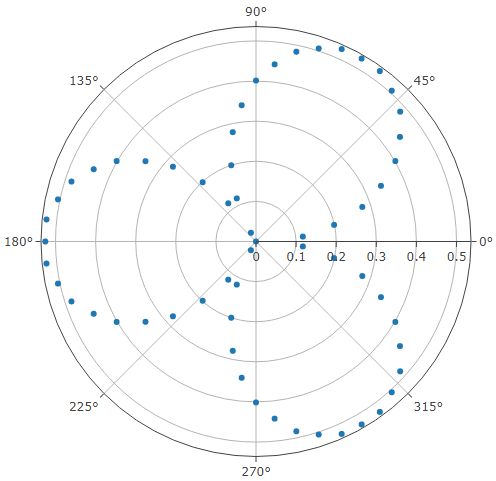}
    \includegraphics[scale = 0.2]{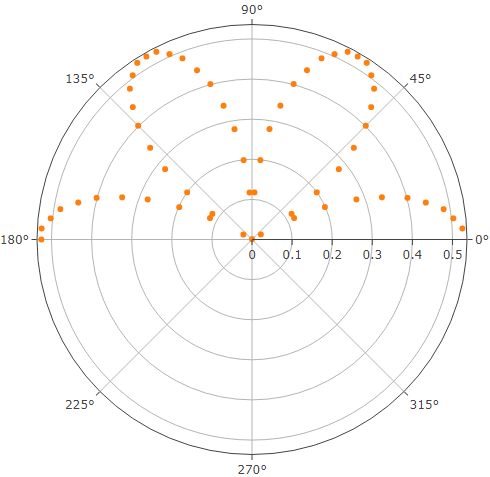}\\
    c) \includegraphics[scale = 0.18]{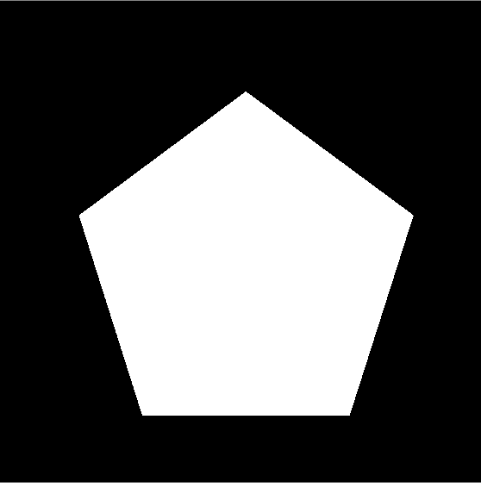}
    \includegraphics[scale = 0.2]{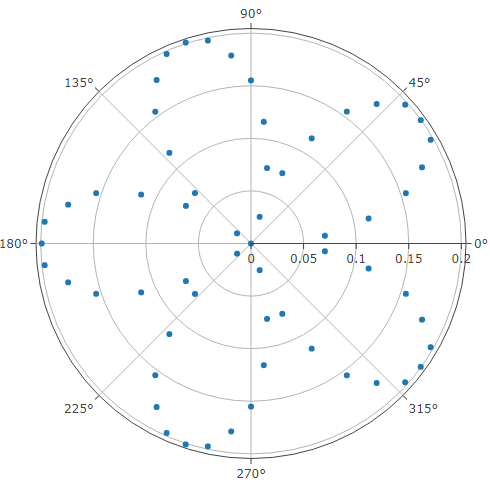}
    \includegraphics[scale = 0.2]{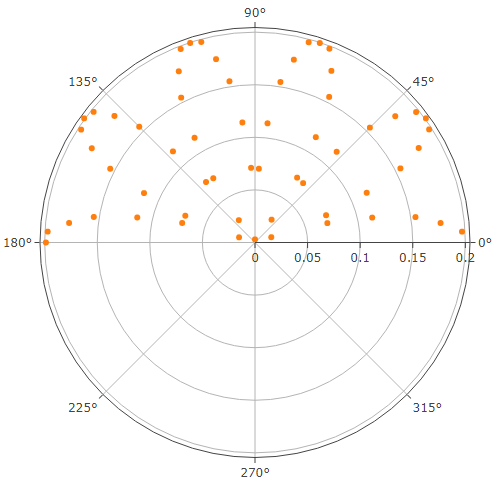}
     d)\includegraphics[scale = 0.3]{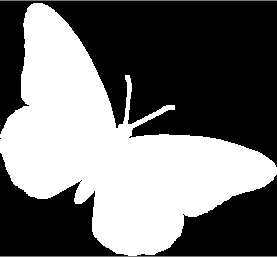}
    \includegraphics[scale = 0.18]{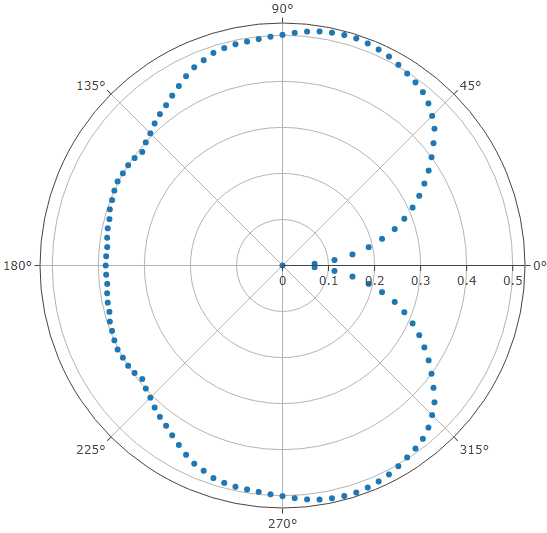}
    \includegraphics[scale = 0.18]{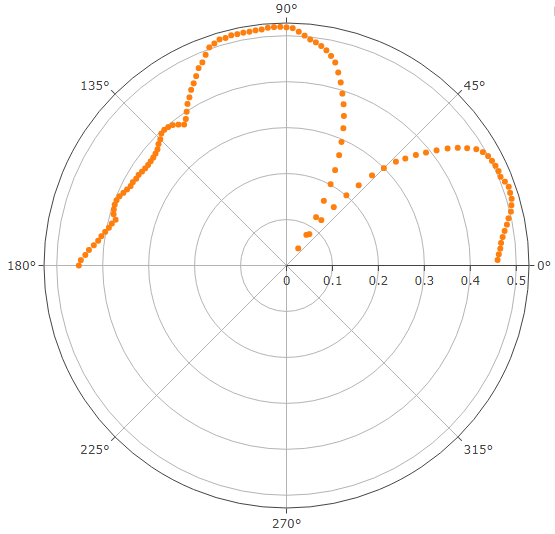}\\
    e)\includegraphics[scale = 0.3]{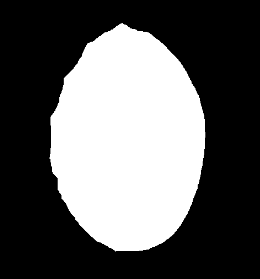}
    \includegraphics[scale = 0.18]{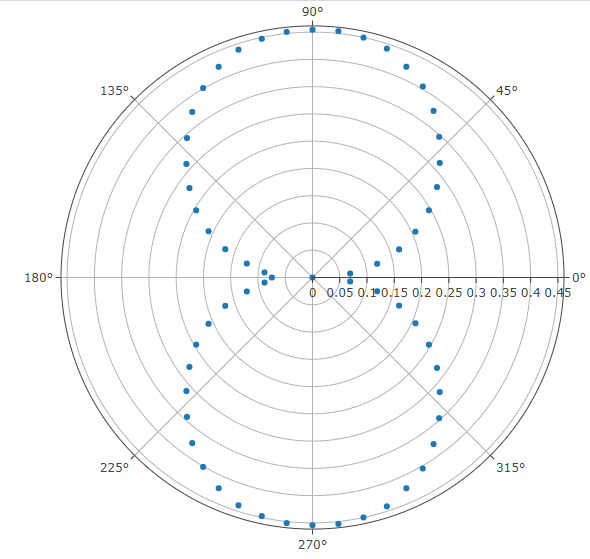}
    \includegraphics[scale = 0.18]{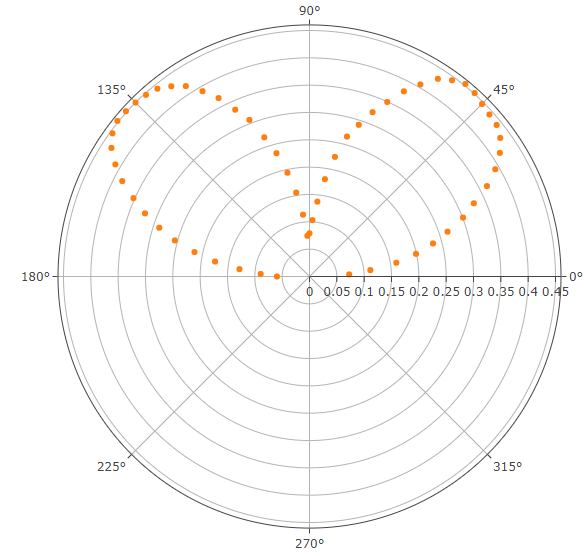}
    f)\includegraphics[scale = 0.22]{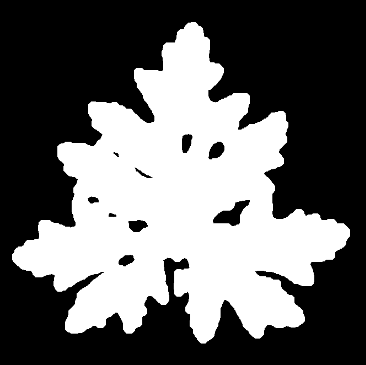}
    \includegraphics[scale = 0.18]{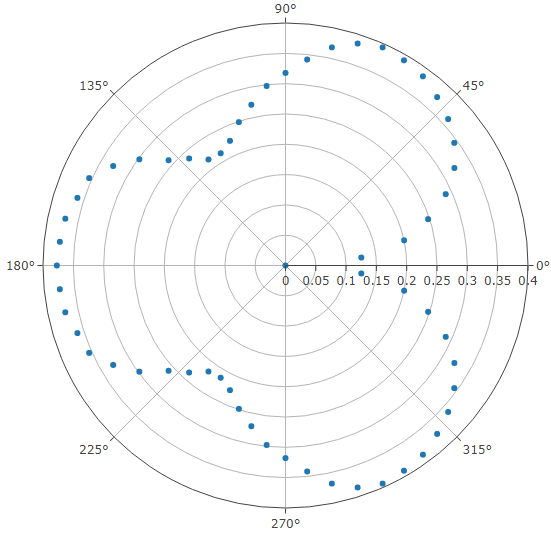}
    \includegraphics[scale = 0.18]{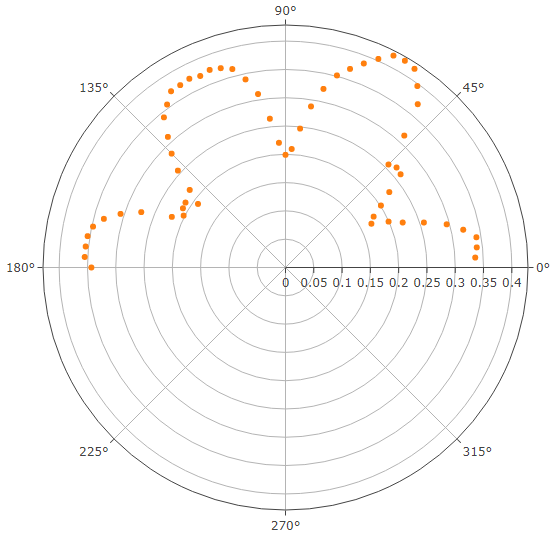}\\
  \caption{Six shapes analysed using the Extended Persistent Homology Transform next to polar plots of the asymmetry score for orthogonal transformations applied to the image. Plots with blue dots at $(\theta,r)$ are for anti-clockwise rotations by angle $\theta$.   
  In the plots with orange dots, the orthogonal transformations studied are reflections in a plane with angle $\theta$ to the positive $x$-axis. 
  The images are a) Circle, b) Triangle, c) Pentagon, d) Butterfly, e) Quercus Suber leaf, f) Geranium species leaf. Note that the radial scale for each shape is different; approximate upper bounds for $r$ are a) 0.055, b) 0.55, c) 0.20, d) 0.55, e) 0.45 f) 0.40. 
  The discretised circle in (a) provides a benchmark tolerance of $0.06$ for accepting a transformation as a symmetry.  Within this tolerance, we see that the triangle and pentagon each have the expected rotation and reflection symmetries: 
  b) has minimal asymmetry scores when the rotation angle is $\theta = 0^\circ, 120^\circ, 240^\circ$  and reflection plane is at $30^\circ, 90^\circ$ and $150^\circ$;   
  c) when the rotation angle is $0^\circ, 72^\circ, 144^\circ, 216^\circ, 288^\circ$ and reflection plane at $18^\circ, 54^\circ, 90^\circ, 126^\circ$ and $162^\circ$. 
  The butterfly in d) has no non-trivial rotational symmetry, but a reflection symmetry in a  plane at angle $\theta = 55.5^\circ$. 
  The leaf in e) has minima of about $0.07$ in rotational asymmetry scores when $\theta = 0^\circ$ and $180^\circ$ and scores of $0.08$ for reflection in planes at $0^\circ, 90--93^\circ, 180^\circ$. 
  In f), we see that none of the orthogonal transformations are near the threshold for a `true' symmetry, but there are still local minima in the scores at rotation angles of $120^\circ, 240^\circ$ and reflections in planes with angles $27^\circ, 90^\circ$ and $144^\circ$, capturing the approximate triangular shape of the leaf.}
  \label{fig:teaser}
}

\abstract{Symmetry is ubiquitous throughout nature and can often give great insights into the formation, structure and stability of objects studied by mathematicians, physicists, chemists and biologists. However, perfect symmetry occurs rarely so quantitative techniques must be developed to identify approximate symmetries. To facilitate the analysis of an independent variable on the symmetry of some object, we would like this quantity to be a smoothly varying real parameter rather than a boolean one.  
The extended persistent homology transform is a recently developed tool which can be used to define a distance between certain kinds of objects. Here, we describe  how the extended persistent homology transform can be used to visualise, detect and quantify certain kinds of symmetry and discuss the effectiveness and limitations of this method. %
} 

\CCScatlist{
  \CCScatTwelve{Human-centered computing}{Visualization}
  {Visualization techniques}{}; \CCScatTwelve{Mathematics of computing}{Continuous mathematics}{Topology}{Algebraic topology}{}
}

\begin{document}


\firstsection{Introduction}

\maketitle

One of the key goals in applied topology is finding meaningful ways to understand and characterise shape. Tools such as persistent homology and the Euler curve transform have been shown to create useful shape descriptors in a variety of applications including bone microstructures \cite{bones}, brain tumours \cite{brain}, and even abstract networks \cite{networks}. 
A less well-known tool from algebraic topology is extended persistence, introduced in~\cite{XPH}. 
This has recently been developed into the extended persistent homology transform (XPHT) in~\cite{XPHT}. 
For suitably nice subsets of Euclidean space, the XPHT quantifies topological changes in height filtrations taken over all directions, and provides a way to measure distances between two shapes. 


The identification of symmetries is important in the natural sciences due to a connection with the stability, structure and assembly of various objects of interest. 
However, most symmetry in biological or chemical contexts is imperfect so methods that can quantify approximate symmetry are required. 
An ideal method for detecting and quantifying approximate symmetry would be applicable to a wide variety of objects, and a broad range of symmetry operations, be robust to noise, and require little pre-processing on the part of the researcher.  

In this paper we demonstrate how the extended persistent homology transform (XPHT) can be used to detect and quantify the symmetry of planar shapes and assess its suitability as a method for detecting and quantifying approximate symmetries.

\subsection{Existing Methods and Related Work}

The problem of understanding approximate symmetry has been well studied and many solutions have been proposed; a number of symmetry detection methods are discussed in the review paper \cite{SIG}. We briefly describe three approaches that are closest to our method.

The Planar Reflective Symmetry Transform (PRST) was developed in \cite{PRST}. 
It quantifies how symmetric an object is under reflection planes, specified by their orientation and location within the image domain, using the $L_2$ norm of the difference between a function defining the object and a symmetric analogue.  
This method is useful as it does not require any guidance from the researcher and can easily work on 2D and 3D models, but is limited by only assessing reflection symmetries. It also has a reported running time of approximately 40 seconds when running on a 64x64x64 voxel grid, which makes it very reasonable for use on smaller 3D datasets.

A flexible method for symmetry quantification described in \cite{TIA} uses Transform Information, an adaption of the Kullback-Liebler Divergence developed in \cite{TI}.  This method can use shape information as well as coloration to detect a broad range of symmetries from photographs of biological samples or simulated data. 
It can be applied to a wide range of applications, but requires  specifications from the researcher to ensure a sensible set of transformations are being assessed.  

Another approach focuses on identifying symmetric sub-parts of an image, such as the work of \cite{2011} and \cite{AEG}. 
Those authors find symmetries in scalar fields  using combinatorial descriptors of the level sets, namely the contour tree in \cite{2011}, and Morse-Smale cell decomposition in \cite{AEG}. 
%
%
These approaches have the advantage of not needing to specify a transformation, 
but may require substantial user interaction, notably in pre-processing the data to stabilise it with respect to input noise in \cite{2011} and setting thresholds when clustering in \cite{AEG}.

The XPHT method we describe in this paper requires little pre-processing, like the PRST, and can assess rotation as well as reflection symmetries of 2D objects. 
Our method does not have the flexibility of Transform Information, but the centering procedure described in section~\ref{remarks} automatically selects the origin for rotation and reflection planes at a location that is most likely to highlight symmetry. 
Some researcher input would be required for our method to detect partial symmetries in a scalar field, as in~\cite{AEG}. For each level set threshold selected, the XPHT procedure can be applied to quantify the presence of rotational and reflection symmetries. 



\subsection{Outline}

This paper has three main parts. In Section~\ref{GuidetoXPHT} we briefly discuss the theory behind the Extended Persistent Homology Transform (XPHT) and how we use the XPHT to measure the difference between an object and its image under Euclidean isometries (specifically, elements of the orthogonal group $O_n(\R)$). 
In Section~\ref{Results}, we present the results of three tests that were undertaken to show the uses of this new methodology for symmetry quantification. The first test shows this method can accurately detect known symmetries using an appropriate symmetry threshold and identify approximate symmetries in more complex images. 
The second test analyses the potential effects of the resolution of the image on the perceived symmetry.
The final test quantifies symmetry in two sets of real-world data: one investigates the bilateral symmetry of different species of leaves and the other studies symmetry in serif and sans-serif fonts. 
Section~\ref{Discussion} discusses the effectiveness and utility of the XPHT symmetry quantification based on the results from the three tests.

\section{A brief guide to the Extended Persistent Homology Transform}
\label{GuidetoXPHT}

The persistent homology transform (PHT), as introduced in~\cite{PHT}, is a powerful method for comparing geometric differences between shapes that have the same essential homology, i.e., the same number of components, handles, voids, etc. For shapes in a particular class (tame, compact, piecewise-linear subsets of $\R^n$), the PHT is injective in the sense that two shapes have the same PHT if and only if they are identical as subsets of $\R^n$. 
If two shapes have different essential homology, the distance between their PHTs becomes infinite, and they can no longer be meaningfully compared. 
The extended persistent homology transform solves this issue and provides a method for comparing shapes with different essential homology. In this section, we briefly define extended persistence diagrams, the extended persistent homology transform (XPHT), a distance between XPHTs of two shapes, and present algorithms for using this information to quantify symmetries of a given object.

\subsection{Extended Persistence}


For the purposes of this paper we use the term `object' to describe a piecewise-linear embedding of a finite simplicial complex in $\R^n$. In regular persistent homology, we might study such an object $M \subset \R^n$ using a height function $h_v:M\to \R$   defined by taking the dot product of points in $M$ with respect to a chosen unit vector $v \in S^{n-1}$ i.e., for $x\in M$, $h_v(x) = x\cdot v$. 
We use this function to define a filtration $\{M_t\}_{t\in \R}$ of $M$ by the sublevel sets $M_t = h_v^{-1}((-\infty,t])$.  
Persistent homology associates an interval of $\R$ with each topological feature seen in the sublevel sets of $M$.  
These intervals are defined using the linear maps on homology groups induced by inclusion $i:M_s\to M_t$ for $s\leq t$. For more details about this process one may consult \cite{PHO} or \cite{PHT}.


One inconvenient feature of this method is that the essential homology classes of our object $M$ will be associated to intervals of the form $[a,\infty)$ as they exist in all $M_t$ for $t$ sufficiently large. 
Having intervals whose endpoints are at $+\infty$ creates problems for metrics defined using this method. 
Extended persistent homology, initially defined in \cite{XPH} for Morse functions on manifolds, resolves this by coupling the sublevel set filtration defined by $h_v$ to a descending filtration, $\{M^t\}_{t\in \R}$, of superlevel sets  $M^t = h_v^{-1}([t,\infty))$. 
The topological structure of the descending filtration is quantified 
using the relative homology of the pair $(M,M^t)$ as discussed in \cite{PHO} and \cite{DPC}. 
The resulting sequence of homology groups is 
\[ H_k(\emptyset,\emptyset) \to  H_k(M_t,\emptyset) \to \cdots 
\to H_k(M, M^t) \to \cdots H_k(M, M).  \]
where $t$ is an increasing parameter in the first half of the sequence and a decreasing parameter in the second.  
This results in finite length intervals being associated to all topological features of $M$. 
An example of an extended persistent homology calculation is depicted in Figure~2.  A more detailed definition of the parameter space is given in the following subsection. 


\begin{figure}
    \centering
    \includegraphics[scale=0.35]{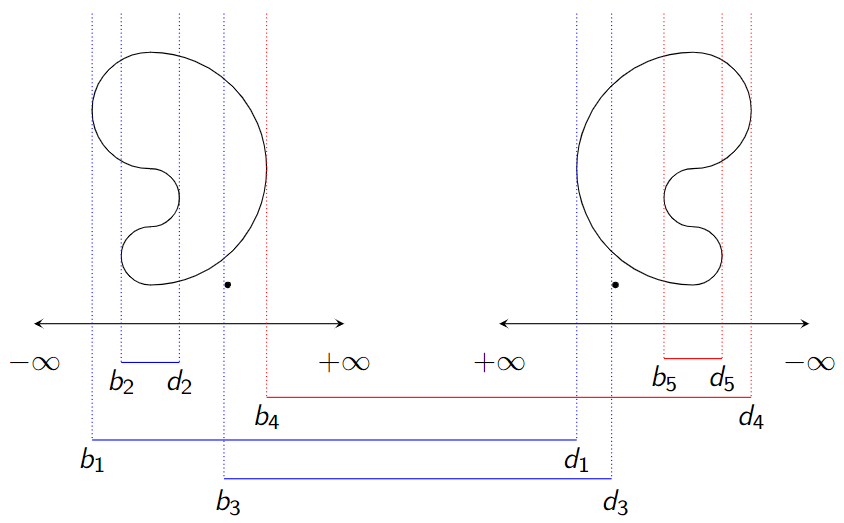}
    \caption{A diagram depicting the extended persistent homology filtration for an embedding of $S^1$ in $\R^2$, together with a small additional component. On the left we see the object $M$ and the filtration parameter $(s, \textbf{Ord})$ along the $x$-axis. The height function $h_{v}$ has $v = (1,0)$ and returns the $x$-coordinate of points in $M$.  
    On the right we duplicate the object $M$ but with the $x$-axis flipped so that the filtration parameter value $(s, \textbf{Rel})$ now decreases going from left to right.  The parameters $b_i$ mark where topological features are ``born'' during our filtration and $d_i$ show where those topological features ``die''. Note that as real numbers, $b_1 = d_4 < b_2 = d_5 < d_2 = b_5 < b_3 < d_3 < b_4$.  The blue intervals denote elements of the $0$th extended persistence diagram and the red intervals denote elements of the $1$st extended persistence diagram. Intervals contained entirely in the ascending filtration on the left hand side are classes that are born and die during ordinary persistent homology. Intervals contained entirely in the descending filtration on the right hand side are classes that are born and die during the relative persistent homology. Essential classes, those in the homology groups of the whole space, have endpoints that cross between the ascending and descending filtration.}
    \label{fig:XPHTCalc}
\end{figure}

The collection of intervals found using degree-$k$ homology groups of the ascending and descending filtration of $M$ by the height function $h_v$ is called the $k$th extended persistence diagram and denoted $XPH_k(M,v)$.

It might seem more straightforward to obtain finite length intervals by simply truncating the regular persistent homology filtration at some maximum value, $H$.  This would mean essential classes of $M$ have an interval of the form $\left[a_i, H\right)$. However, this makes every essential class take the same death value, irrespective of the spatial extent of that topological feature.  Using the extended persistent homology sequence, each essential class is given a death value corresponding to the appropriate level set threshold. Small noise-induced holes or additional components will have suitably short intervals, and have less impact on distance calculations.  

For the example in Figure~\ref{fig:XPHTCalc}, a truncation approach would yield intervals $\left[b_1,H\right)$, $\left[b_2,d_2\right)$, $\left[b_3,H\right)$ in degree-0, and $\left[b_4,H\right)$ in degree-1. In comparison, the $0$th extended persistence diagram has ordinary interval $\left[b_2,d_2\right)$, and essential intervals $\left[b_1,d_1\right)$ and $\left[b_3,d_3\right)$ while $XPH_1(M,x)$ has a relative interval $\left[b_5,d_5\right)$ and essential interval $\left[b_4,d_4\right)$. 
Since we must have $H \geq b_4=d_1$, we see that the $0$-th essential class associated with the large bean may not change too much, but the interval for the $1$-cycle could become vanishingly small, and the interval corresponding to the small spot increases significantly.

\subsection{Distances on the Space of Extended Persistence Diagrams}

By defining a metric for extended persistence diagrams, we obtain a distance between objects that can then be used to quantify their geometric and topological differences. 
The metric we use is an optimal transport one, commonly called a Wasserstein distance in the persistent homology literature. 
In this paper, for the purpose of simplicity, 
we restrict discussion to the 1-Wasserstein distance between extended persistence diagrams, denoted $W_1(X,Y)$; see~\cite{XPHT} for a more general formulation. The $1$-Wasserstein distance is chosen because, compared to other $p$-Wasserstein distances, it is more closely related to the area (or volume for 3D shapes) of the symmetric difference between $X$ and $Y$.

To understand this metric, we must first describe the intervals formed by the extended filtration. We define the parameter space $\Theta$ using a pair $(r,\textbf{Type})$, where $r$ is a real number and $\textbf{Type} \in \{\textbf{Ord},\textbf{Rel}\}$. Points of the form $(s,\textbf{Ord})$ refer to the ascending filtration of sublevel sets and points of the form $(s,\textbf{Rel})$ to the descending filtration where we apply relative homology.
We define a total order, $\leq$, on $\Theta$ using the following rules:
\begin{align*}
    (s,\textbf{Ord}) \leq (t,\textbf{Ord}) &~~\text{if}~~ s\leq t \in (\R,\leq)\\
    (s,\textbf{Rel}) \leq (t,\textbf{Rel}) &~~\text{if}~~ t\leq s \in (\R,\leq)\\
    (s,\textbf{Ord}) \leq (t,\textbf{Rel}) &~~\text{for all}~~ s,t \in \R.
\end{align*}
The intervals in extended persistence diagrams will always be of the form $[\mathfrak{b},\mathfrak{d})$ for $\mathfrak{b},\mathfrak{d}\in \Theta$ such that $\mathfrak{b} \leq \mathfrak{d}$. 

We next define a distance $d_\Theta:\Theta^2 \to \R\cup \{\infty\}$ between elements of $\Theta$ using the following rules:
\begin{align*}
    d_\Theta((s,\textbf{Ord}),(t,\textbf{Ord})) &= |s-t|\\
    d_\Theta((s,\textbf{Rel}),(t,\textbf{Rel})) &= |s-t|\\
    d_\Theta((s,\textbf{Ord}),(t,\textbf{Rel})) &= \infty
\end{align*}
A distance, $d_{\mathfrak{I}}$, between intervals simply adds the distance between corresponding endpoints:
\[d_{\mathfrak{I}}\left([\mathfrak{b}_1,\mathfrak{d}_1),[\mathfrak{b}_2,\mathfrak{d}_2)\right) = d_\Theta(\mathfrak{b}_1,\mathfrak{b}_2) + d_\Theta(\mathfrak{d}_1,\mathfrak{d}_2).\]

Before we define the metric on diagrams, we must first introduce a set of `ephemeral intervals', $\textbf{Eph}$, defined as
\begin{align*}
    \textbf{Eph} :=& \{[(s,\textbf{Ord}),(s,\textbf{Ord})) ~|~ s\in \R\} ~\cup \\
    &\{[(s,\textbf{Rel}),(s,\textbf{Rel})) ~|~ s\in \R\} ~\cup\\
    &\{[(s,\textbf{Ord}),(s,\textbf{Rel})) ~|~ s\in \R\}.
\end{align*}

We will notate an extended persistence diagram as a collection of intervals $X = \{[\mathfrak{b}_s,\mathfrak{d}_s) ~|~ s \in S^X\}$ where $S^X$ is some finite index set. 

Now, given two extended persistence diagrams, $X,Y$ we define a transportation plan, $T$, between $X$ and $Y$ to be a triple $(\hat{S}^X, \hat{S}^Y,\rho)$ where $\hat{S}^X \subset S^X$, $\hat{S}^Y \subset S^Y$ and $\rho:\hat{S}^X \to \hat{S}^Y$ is a bijection. We define the cost of a transportation plan $c(T)$ as follows.
\begin{align*}
    c(T) &= \sum_{s \in \hat{S}^X} d_{\mathfrak{I}}([\mathfrak{b}_s,\mathfrak{d}_s),[\mathfrak{b}_{\rho(s)},\mathfrak{d}_{\rho(s)}))\\
    &+ \sum_{s \in S^X\setminus \hat{S}^X} \inf_{I \in \textbf{Eph}} \{d_{\mathfrak{I}}([\mathfrak{b}_s,\mathfrak{d}_s),I)\}\\
    &+ \sum_{s \in S^Y\setminus \hat{S}^Y} \inf_{I \in \textbf{Eph}} \{d_{\mathfrak{I}}([\mathfrak{b}_s,\mathfrak{d}_s),I)\}
\end{align*}
With all of this established, we can finally define the 1-Wasserstein metric on the space of extended persistence diagrams. Given extended persistence diagrams, $X,Y$ we set \[W_1(X,Y) := \inf\{c(T) ~|~ ~\text{$T$ is a transportation plan between $X$ and $Y$}\}.\]

Further details and the generalisation to $p$-Wasserstein distances and other metrics on extended persistence diagrams can be found in \cite{XPHT}. 

\subsection{What is the XPHT?}

The Extended Persistent Homology Transform (XPHT) maps the space of objects $M \subset \R^n$ to parameterised $n$-tuples of extended persistence diagrams.  The parameter here is a unit vector $v \in S^{n-1}$, and extended persistence diagrams are computed from the height filtrations $h_v$ defined earlier.  We write 
\[ XPHT(M,v) = (XPH_0(M,v),...,XPH_{n-1}(M,v)). \]


We use this transform to define a notion of distance between objects, $M, N \subset \R^n$, first by extending the 1-Wasserstein distance between diagrams to $n$-tuples of diagrams, then integrating over all possible direction vectors: 
\[ W_1(XPHT(M,v),XPHT(N,v)) = \sum_{i=0}^{n-1} W_1(XPH_i(M,v),XPH_i(N,v)) \]
and 
\[d^{\text{XPHT}}(M,N) = \int_{S^{n-1}} W_1(XPHT(M,v),XPHT(N,v)) dv.\]


This definition of distance is technically a psuedo-metric rather than a metric, but it is very close to being a true metric. If two shapes have distance zero then they differ only by some essential classes that lie entirely within level sets regardless of direction. In particular, for compact $M, N$, if $d^{\text{XPHT}}(M, N)=0$ then the symmetric difference between $M$ and $N$ is a set of isolated points. 
If we restrict our analysis to subsets constructed from pixelated images as described in~\cite{XPHT} then $d^{\text{XPHT}}$ is a metric because each pixel has a positive width.

\subsection{How does this allow us to quantify asymmetry?}

Now that we have a  distance between objects we can use this to compare an object, $M$ to its image under a transformation $T$, creating what we call the asymmetry score of $M$ under $T$, 
\[S(M,T) = d^{XPHT}(M,T(M)).\]
Calculating this score involves first calculating the XPHT for $M$ and $T(M)$, however, in the situation where $T \in O_n(\R)$ we can calculate this more efficiently by using the fact that $T$ is a transformation of $M$ as well as an isometry of $S^{n-1}$.

To see how this works, we first note that \[XPH_k(M,v) = XPH_k(T(M),T(v)).\] This is because the sublevel sets, $M_t$, are exactly the points $x\in M$ such that $\langle x, v\rangle \leq t$, which also satisfy $\langle (T^{-1}\circ T) (x), v\rangle \leq t$ as $T^{-1}\circ T = I$. Hence, from the fact that the transpose of $T$ is its inverse, these points must also satisfy $\langle T(x), T(v) \rangle \leq t$ which is exactly the sublevel set, $T(M)_t$ with respect to the function $h_{T(v)}$.
We use this result to replace $XPHT(T(M),v)$ with $XPHT(M,T^{-1}v)$ and rewrite our asymmetry score as
\[S(M,T) = \int_{S^{n-1}} W_1(XPHT(M,v),XPHT(M,T^{-1}(v))) ~dv .\]
However, it is preferable for this to be written in terms of $T$ rather than $T^{-1}$. We note that $O_n(\R)$ is the isometry group of $S^{n-1}$ and the  integral is taken over all of $S^{n-1}$, so we can make a change of variable $u=T^{-1}(v)$ without changing the value of the integral. Using this and the symmetry of $W_1$ we obtain:
\[S(M,T) = \int_{S^{n-1}} W_1(XPHT(M,u),XPHT(M,T(u))) ~du .\]

This score has the property that $S(M,T) = 0$ when $T$ is a symmetry of $M$ because if $T(M)=M$ then $S(M,T) = d^{XPHT}(M,M) = 0$. This suggests that the further the score is away from $0$ the more asymmetric the object is under the given transformation.

A discretised version of this asymmetry score for 2D objects is given by \[S(M,T) = \frac{1}{N}\sum_{i=1}^N W_1\left(XPHT\left(M,v_i\right),XPHT\left(M,Tv_i \right)\right).\] 
Here, $N$ is a user-chosen number of directions and \[v_i = \begin{bmatrix} \cos(\frac{2\pi i}{N})\\ \sin(\frac{2\pi i}{N})\end{bmatrix} ~\text{for}~ i \in \{1,...,N\}.\] 
An R-package that computes the XPHT and Wasserstein distances from 2D binary images is available on github~\cite{R}. 
This implementation requires $N$ to be even. 

The symmetry quantification procedure calculates the discretised asymmetry scores for all transformations $T \in O_2(\R)$ that map the set of direction vectors, $\{v_i ~|~ i \in \{1,...,N\}\}$, onto itself.
These are rotations by $\frac{2\pi i}{N}$ for $i \in \{1,...,N\}$ and reflections about the planes that make an angle of $\frac{\pi i}{N}$ with respect to the positive direction of the $x$-axis. 
Our algorithm starts by calculating $XPH(M,v_i)$ for each of the specified direction vectors then computes the pairwise distance matrix $D = [\delta_{ij}]$ using the $1$-Wasserstein distances between the different directions so that \[\delta_{ij} = W_1\left(XPHT\left(M,v_i\right),XPHT\left(M,v_j\right)\right).\]


Algorithm~1 shows how we determine the asymmetry scores associated to the $N$ rotations from this distance matrix.

\begin{algorithm}
\caption{Symmetry Score for Rotations}
\begin{algorithmic} 
\REQUIRE Distance matrix, $D[i,j]$, and the number of equally spaced directions, $N$, used to calculate the XPHT.
\ENSURE 
\STATE SCORES $\leftarrow 1\times N$ matrix
\STATE ANGLES $\leftarrow 1\times N$ matrix 
\FOR{$i = 1:N$}
\STATE SCORES$[i]  \leftarrow 0$
\STATE ANGLES$[i] \leftarrow \frac{360}{N}(i-1)$
\FOR{$j=1:N$}
\STATE SCORES$[i] \leftarrow$ SCORES$[i] + D[j,(i+j-2)~\textbf{mod}~N + 1]$
\ENDFOR
\ENDFOR
\RETURN $2\times N$ Matrix $\begin{bmatrix} \text{SCORES}[i]\\\text{ANGLES}[i] \end{bmatrix}$ 
\end{algorithmic}
\end{algorithm}

This works because $i$ sets the angle between two initial direction vectors, and incrementing $j$ accumulates the distances between diagrams with this same angle between direction vectors. The -$2$ in our iterative step is to account for the fact that arrays in R are indexed starting at $1$ not $0$.


\begin{algorithm}
\caption{Symmetry Score for Reflections}
\begin{algorithmic} 
\REQUIRE Distance matrix, $D[i,j]$, and the number of equally spaced directions, $N$, that the XPHT was calculated for.
\ENSURE 
\STATE SCORES $\leftarrow 1\times N$ matrix
\STATE ANGLES $\leftarrow 1\times N$ matrix 
\FOR{$i = 1:N$}
\STATE SCORES$[i]  \leftarrow 0$
\STATE ANGLES$[i] \leftarrow \frac{180}{N}i$
\FOR{$j = 1:N$}
\STATE SCORES$[i] \leftarrow$ SCORES$[i] + D[j,(i-j-1~ \textbf{mod} N) + 1]$
\ENDFOR
\ENDFOR
\RETURN $2\times N$ Matrix $\begin{bmatrix} \text{SCORES}[i]\\\text{ANGLES}[i] \end{bmatrix}$ 
\end{algorithmic}
\end{algorithm}

Algorithm~2 shows how we collect the asymmetry scores associated to the $N$ possible reflection transformations.

The algorithm works because fixing a position in our distance matrix fixes two direction vectors, for which there is a unique reflection plane that maps these vectors onto each other. If we rotate one vector clockwise by a certain angle and the other vector counterclockwise by the same angle, then the reflection plane that maps our rotated direction vectors onto each other is the same as the one that mapped our initial vectors onto each other. This is exactly what is achieved by moving one entry down and one entry to the left in our distance matrix, using periodic boundaries. Our offset i fixes a particular point in the distance matrix, defining a reflection plane, and then as j increases we sum over all entries whose defining vectors are mapped to each other by the same reflection plane.

The algorithms have been set up so that the rotation scores in our final matrix are positioned above the angle that defines each rotation and the reflection scores are above the angle the plane that defines the reflection makes with the positive direction of the $x$-axis. We can then plot these outputs directly, or convert to polar coordinates as in Figure~\ref{fig:teaser}.

\subsection{Computational complexity}

For $M\subset \R^2$ given in the form of a binary image the computation time of $XPHT(M,v_i)$ is $O(n \alpha(n))$ where $n$ is the number of pixels on the border between the $0$- and $1$-labelled pixels and $\alpha(n)$ is the inverse of the Ackermann function which grows so slowly it is negligible. Details for this are in \cite{XPHT}. As $XPHT(M,v_i)$ will need to be computed in each of the $N$ directions separately, this implies the total computation of the $XPHT(M)$ is $O(n\alpha(n)N)$.

There are a variety of methods to compute the $1$-Wasserstein distance between two persistence modules $X$ and $Y$. One method is to use Munkres algorithm to find an optimal transportation plan. This algorithm is $O(m^3)$ where $m$ is the maximum number of off-diagonal points in $X$ and $Y$. 

We are interested in the pairwise distances for each pair of directions $\{(v_i,v_j)\mid 1\leq i,j \leq N\}$. As these need to be computed separately this totals $O(N^2m^3)$ to compute the matrix $D$.  From these pairwise distances we can compute the symmetry score matrices for rotation and reflection. Each of these matrices takes $O(N^2)$ time to compute from the entries of $D$.   

\subsection{Remarks concerning the reflections and rotations considered.}
\label{remarks}

Transformations $T \in O_n(\R)$ fix the origin; its location with respect to an object will have a significant effect on whether $T(M) = M$. Here we describe how this origin is set by finding a central point using the algorithm described in \cite{PHT}.
This method acts directly on the persistence diagrams, using the height of the first birth in $XPH_0(M,v_i)$ for each direction $v_i$ to identify points on the outer boundary(ies) of the object. A normalised vector sum of these points defines the center of the object. This is conceptually similar to determining the centre of a regular polygon as the mean of the sum of its vertices.

\begin{figure}
    \centering
    \includegraphics[scale=0.4]{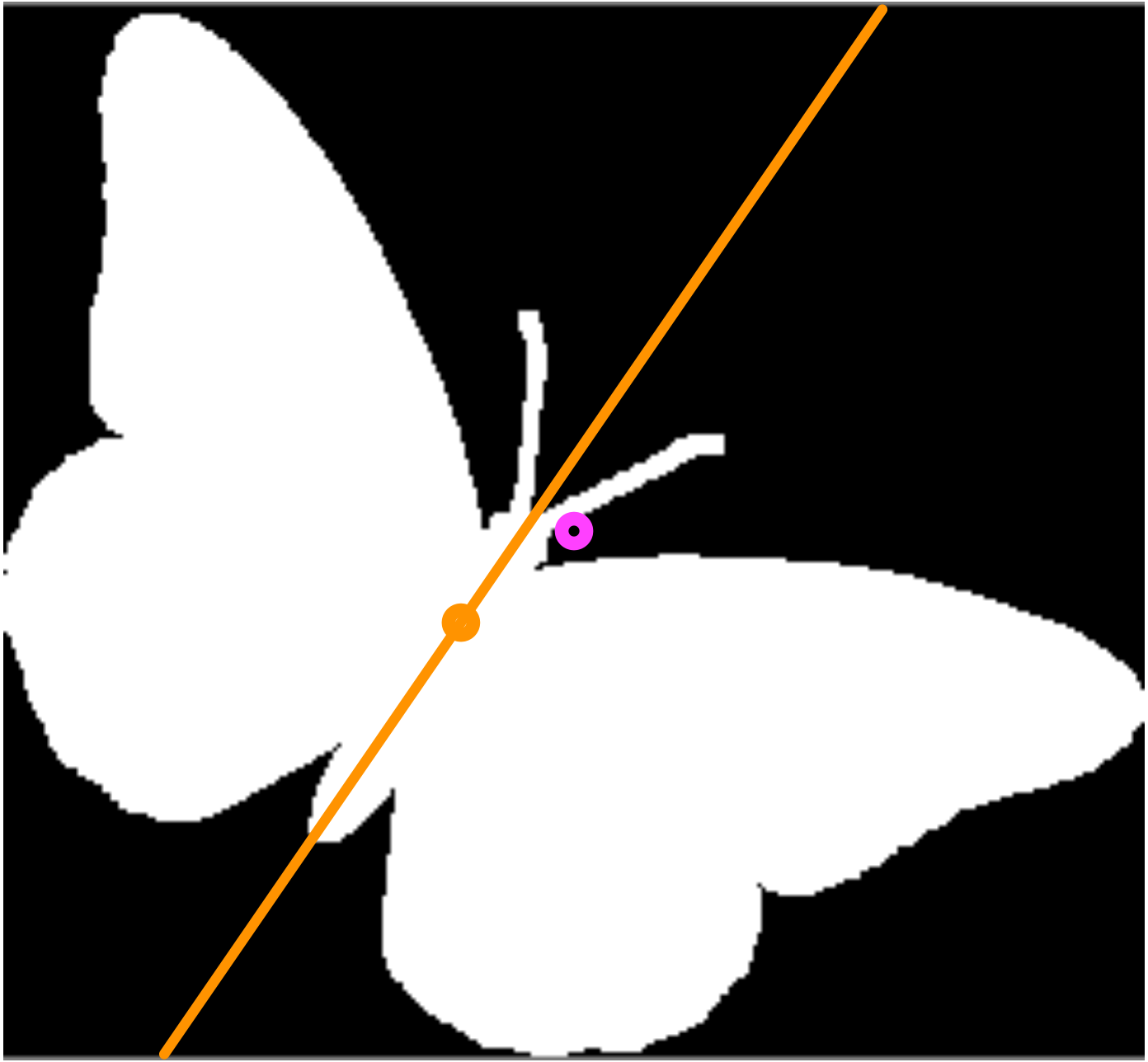}
    \caption{The image from Figure 1.d with the center of the image domain marked with a magenta `o' and the derived center marked in orange. The reflection plane, also drawn in orange, passes through the derived center at $55.5^\circ$, the angle for which the asymmetry score is minimal.}
    \label{fig:Butterflycenter}
\end{figure}

The reflection transformations we construct are in planes that pass through the derived centre and the rotations also fix that point. Of course the derived centre will not necessarily be at the centre of the image domain, nor at the origin of a predefined coordinate system. An illustrative example is shown in Figure~\ref{fig:Butterflycenter}.

While this limits the number of possible symmetry actions we consider with our method, it also focuses our attention on the transformations that are likely to be symmetries of the object. 
This is because when an object is symmetric under some rotation, the points on the convex hull of the object should map to other points on the convex hull at the same distance from the centre of rotation. A similar argument holds for reflections. 


We also note that our method cannot be used directly to detect partial symmetries of an image, where subparts of the image are compared to each other, as the algorithm works with the object as a whole. There are ways to isolate parts of an image and compare them using the XPHT but this will likely require significant preprocessing of data on the part of the researcher.

\section{Results}
\label{Results}

Having developed a program to implement the extended persistent homology transform on binary images and associated code to compute the asymmetry score of a number of orthogonal transformations of these images, we perform a number of tests to determine how appropriate it is for practical uses. 

\subsection{Visualising known symmetry}

The first test is to see if our methodology can detect and quantify the asymmetry of objects and how it could be used as a way to visualise, quantify and compare the asymmetry of certain objects under a collection of transformations. For this we sourced images with known symmetries from the MPEG-7 Shape Matching Dataset \cite{MPEG} (Figure 1.a-d), and leaf images exhibiting approximate symmetry sourced from \cite{leaf} (Figure 1.e-f). We ran these images in our program and have plotted the asymmetry scores of a number of orthogonal transformations in polar plots. We note that part of our algorithm involves rescaling the object in the image so that it just fits inside a unit disc, allowing us to compare the caclulated asymmetry scores of our objects more fairly. The findings are summarised in Figure 1.

\subsection{Resolution Analysis}

\begin{figure}
    \centering
    \includegraphics[scale=0.4]{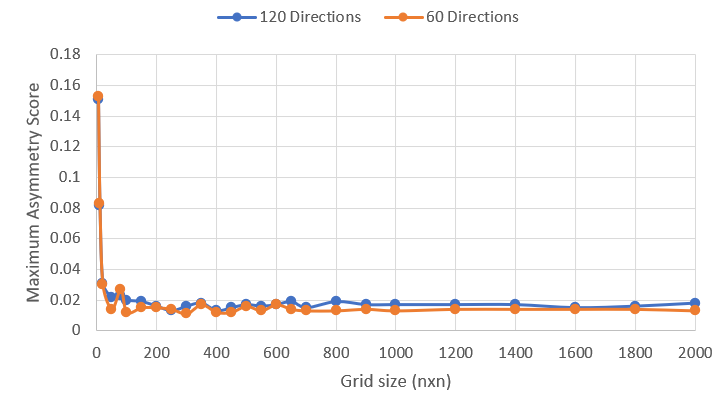}
    \caption{The Maximum Asymmetry Score obtained using our methodology on a range of discretised discs with 60 directions (in orange) and 120 directions (in blue). We see that the Maximum asymmetry score obtained seems fairly stable after the resolution goes past a 100x100 pixel grid.}
    \label{fig:MaxAsymmetry}
\end{figure}

\begin{figure}
    \centering
    \includegraphics[scale=0.5]{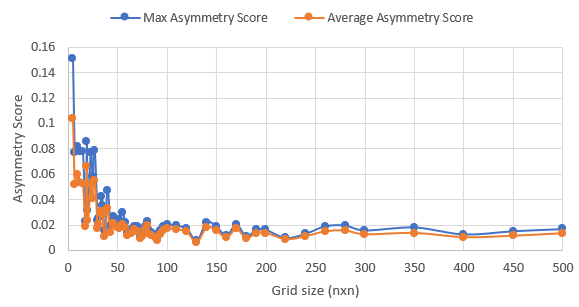}
    \caption{A comparison of the maximum asymmetry score and the average asymmetry score of  discretised discs under the set of orthogonal transformations that leave the set of 120 directions invariant.}
    \label{fig:MaxvAverage}
\end{figure}

\begin{figure}
    \centering
    \includegraphics[scale=0.4]{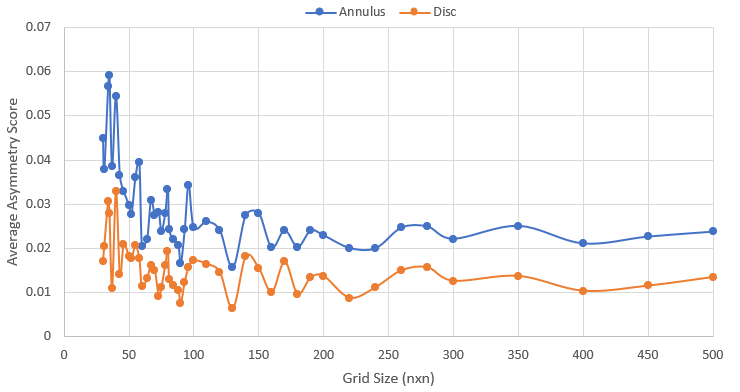}
    \caption{The average asymmetry score of discretised discs and annuluses (circles) under the set of orthogonal transformations that leave our set of 120 directions invariant. We see that the annuluses have average scores that are approximately double those of the discs.}
    \label{fig:DiscvAnnulus}
\end{figure}

Our second test is to see how sensitive the program is to the image resolution. To do this we generate a series of images of discretised discs whose radius is approximately one third the size of the image. 
Using the same set of direction vectors for each image, we calculate the asymmetry values associated to reflection and rotation transformations, for $N=60$ and $N=120$. The plot of maximum asymmetry score across all transformations against the resolution of the image is presented in Figure~\ref{fig:MaxAsymmetry}. 
We also compare the average asymmetry score with the maximum score  using N=120 directions; the results are presented in Figure~\ref{fig:MaxvAverage}. 
Finally, we compare a discretised circle and annulus, objects with a different number of essential homology classes but otherwise similar shape and symmetry.   The results shown in Figure~\ref{fig:DiscvAnnulus} demonstrate that the observed average asymmetry scores are proportional to the number of homology classes. 

\subsection{Bilateral Symmetry of Leaves}

A third test of our methodology is to see whether we can use it to answer some simple research questions. 
In this case we seek to determine which species of leaf exhibits the strongest bilateral symmetry about the mid-vein of the leaf. 

\begin{figure}
    \centering
    \includegraphics[scale=0.45]{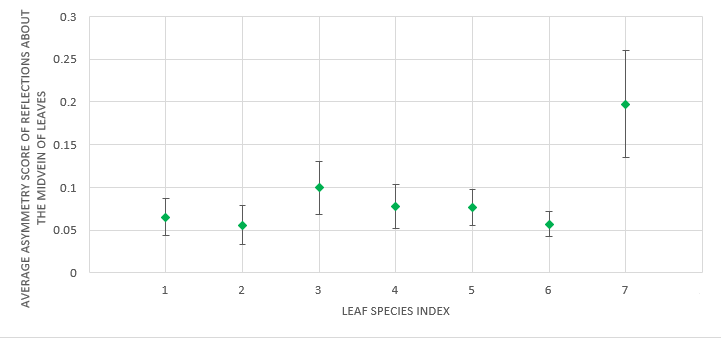}
    \caption{Average asymmetry score obtained for reflection along the midvein for the seven different species of leaf illustrated in Figure 8.}
    \label{fig:LeafPlot}
\end{figure}

For this analysis we use ten images from each of seven different species of leaf, sourced from \cite{leaf}, and run our algorithm using 120 directions. We take the reflection symmetry with the lowest asymmetry score that is approximately in line with the midvein of the leaf as the bilateral symmetry score of the leaf. Lower values imply stronger bilateral symmetry.  The mean values for each leaf species are plotted in Figure~\ref{fig:LeafPlot}. 

\begin{figure}
    \centering
    1. \includegraphics[scale = 0.4]{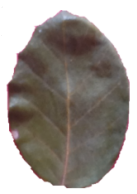}
    2.\includegraphics[scale = 0.4]{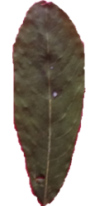} 
    3.\includegraphics[scale = 0.4]{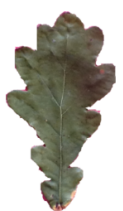} 
    4.\includegraphics[scale = 0.4]{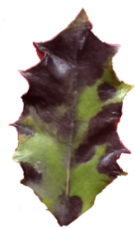}\\ 
    5.\includegraphics[scale = 0.4]{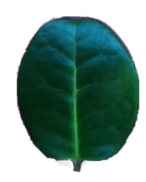} 
    6.\includegraphics[scale = 0.4]{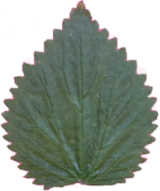} 
    7.\includegraphics[scale = 0.4]{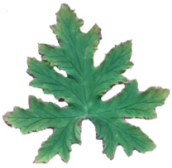}\\
    \caption{1.Quercus suber 2. Salix atrocinera 3. Quercus robur 4. Ilex aquifolium 5. Ilex perado ssp. azorica 6. Urtica dioica 7. Geranium sp.}
    \label{fig:Leaves}
\end{figure}

The leaf species index is 1. \textit{Quercus suber}, 2. \textit{Salix atrocinera}, 3. \textit{Quercus robur}, 4. \textit{Ilex aquifolium}, 5. \textit{Ilex perado ssp. azorica}, 6. \textit{Urtica dioica} and 7. \textit{Geranium sp.}. Example images of these species are depicted in Figure~\ref{fig:Leaves}.

\subsection{Asymmetry of Serif/Sans-Serif Fonts}

The final test of our symmetry quantification is to assess whether serif or sans-serif fonts tend to be more symmetrical.

\begin{figure}
    \centering
    \includegraphics[width=\columnwidth]{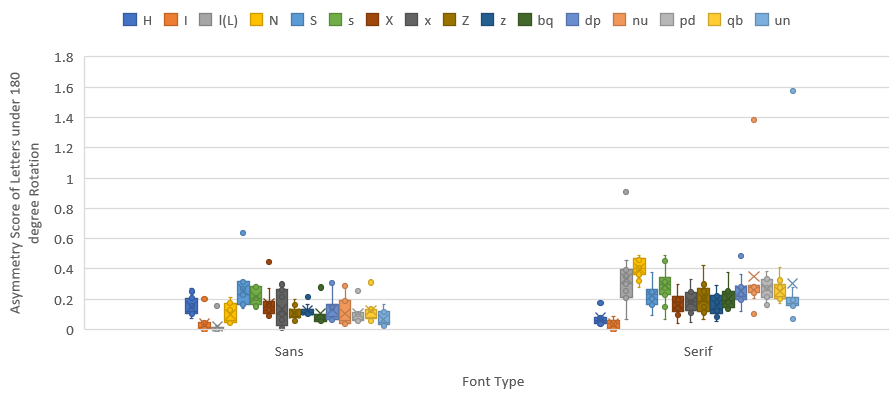}
    \caption{Box and whisker plot of the distribution of asymmetry scores of letters in serif and sans-serif fonts associated to a 180 Degree rotation.}
    \label{fig:Rot}
\end{figure}

\begin{figure}
    \centering
    \includegraphics[width=\columnwidth]{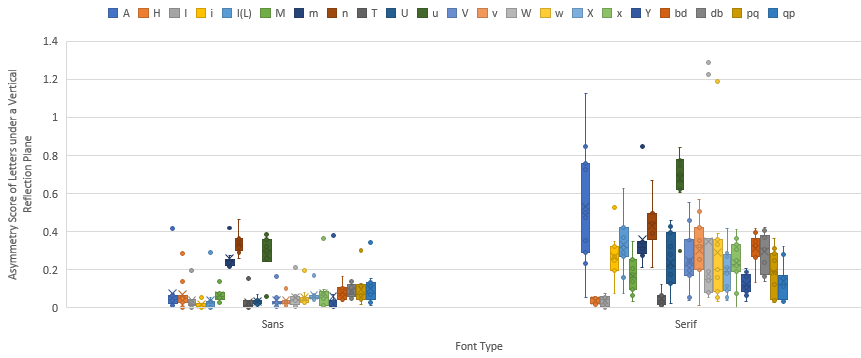}
    \caption{Box and whisker plot of the distribution of asymmetry scores of letters in serif and sans-serif fonts associated to reflections about a vertical axis.}
    \label{fig:vertref}
\end{figure}

\begin{figure}
    \centering
    \includegraphics[width=\columnwidth]{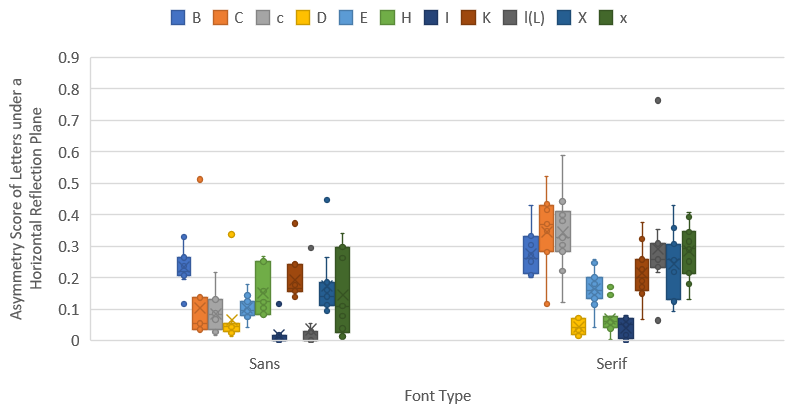}
    \caption{Box and whisker plot of the distribution of asymmetry scores of letters in serif and sans-serif fonts associated to to reflections about a horizontal axis.}
    \label{fig:horizref}
\end{figure}

\begin{figure}
    \centering
    \includegraphics[width=0.8\columnwidth]{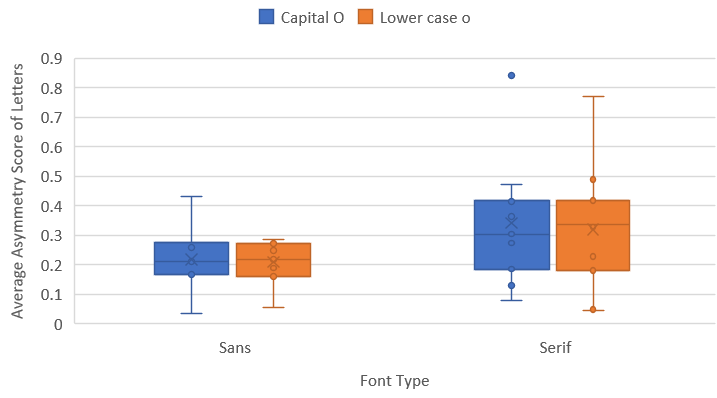}
    \caption{Box and whisker plot of the distribution of the average asymmetry score of letters in serif and sans-serif fonts over all orthogonal transformations that leave a set of 120 direction vectors invariant.}
    \label{fig:lettero}
\end{figure}

For this we begin by choosing 11 serif fonts (Bernard MT Condensed, Bodoni MT, Centaur, Century Schoolbook, Courier New, Elephant, Georgia, Modern No. 20, Perpetua, Rockwell, Times New Roman) and 11 sans-serif fonts (Arial, Bahnschrift, Calibri, Candara, Century Gothic, Comic Sans MS, Corbel, Ebrima, Gadugi, Tahoma, Verdana). We then process images of letters, or pairs of letters, that we expect to be symmetric under 180 degree rotations (H, I, l, N, S, s, X, x, Z, z, bq, dp, nu, pd, qb, un), vertical reflection planes (A, H, I, i, l, M, m, n, T, U, u, V, v, W, w, X, x, Y, bd, db, pq, qp), and horizontal reflection planes (B, C, c, D, E, H, I, K, l, X, x).  
We also quantify the average asymmetry of a capital and lower case `O' under all transformations that leave our set of 120 direction vectors invariant. The results of the 180 degree rotation analysis are presented in Figure~\ref{fig:Rot}, the vertical reflection plane analysis in Figure~\ref{fig:vertref}, the horizontal reflection plane analysis in Figure~\ref{fig:horizref} and the analysis of the letter `o' in Figure~\ref{fig:lettero}.

\section{Discussion}
\label{Discussion}

We shall now consider the results of each of these tests and what they can tell us about the appropriateness and effectiveness of our methodology.

\subsection{Symmetry Visualisation}

From the asymmetry score plots of test images in Figure~\ref{fig:teaser} we see that our methodology works as expected. Of interest here is the magnitude of the asymmetry score in each of the plots as this allows us to compare objects more accurately. 
Note that the objects in these images are scaled to fit within a unit disc, allowing us to directly compare asymmetry scores between these images.
The circular disc is a useful starting point because in the continuous setting they are symmetric under all orthogonal transformations.  As shown in Figure 1.a) all the calculated asymmetry scores for the discretized disc are below a threshold of 0.06.  All other expected symmetries for the simple  images of Figure 1.b-d) also fall below this threshold indicating that our methodology can detect symmetry when present. In Figure 1.e-f) we see how our method can pick up approximate symmetries by looking at local minima of asymmetry scores. Even when the scores do not fall below the 0.06 threshold, we can still see the transformations which come the closest to being symmetries of these shapes.

\subsection{Resolution Analysis}

When analysing the effect of image resolution we discovered a few important features about our method. In Figures~\ref{fig:MaxAsymmetry} and~\ref{fig:MaxvAverage} we see that after reaching a resolution of approximately 100x100 pixels the asymmetry scores of our discretised discs remained stable and approximately the same magnitude regardless of the number of directions used. We conclude that the effects of image resolution are fairly inconsequential. 

In Figure~\ref{fig:MaxvAverage} we investigate our choice to look at a maximum bound on the asymmetry scores of the transformations of our discs by comparing them to the average of the asymmetry scores and note that the averages depict approximately the same trend, though naturally they are slightly lower values than the maxima.

Finally, in Figure~\ref{fig:DiscvAnnulus} we see the effect of increasing the number of essential classes of our objects, from which it seems that in doubling the number of essential classes we have approximately doubled the asymmetry score. This is likely an artifact of the choice to use a 1-Wasserstein distance, analogous to an $L_1$ norm, and that the intervals associated to our essential classes in all directions are approximately the same size so their contribution to the asymmetry score should be approximately equal. Importantly this means that we should account for the number of intervals in an extended persistence diagram when setting a threshold for asymmetry scores. 

\subsection{Symmetry of leaves and letters}

The answers to our simple research questions are easy to glean from the figures we obtained.

The asymmetry scores for leaves show that the species with strongest bilateral symmetry with respect to the midvein is \textit{Salix atrocinera} (2) and the least symmetric species was \textit{Geranium sp.} (7); see Figures~\ref{fig:LeafPlot} and~\ref{fig:Leaves}. The large standard deviation for the Geranium samples is likely due to those leaves having deep lobes and overlapping tips, making the 2D images appear to have holes. As we saw in Figure~\ref{fig:DiscvAnnulus}, a different number of essential classes can significantly impact asymmetry scores. 


From our font analysis we can see that there is some evidence that serif fonts are generally more asymmetric than sans-serif fonts, and this is most strongly observed when considering reflections about a vertical plane.

Perhaps more interestingly, by studying outliers in our font analysis, we can uncover one of the key issues in our methodology. Consider the outlier in `un' and `nu'. The font responsible for this is Modern No.~20 and the feature that distinguishes it from other fonts is that the serifs create a closed loop in the letter n, as seen in Figure 13.

\begin{figure}
    \centering
    \includegraphics[scale=0.5]{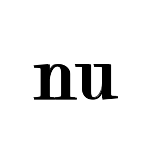}
    \caption{`nu' in the font Modern No. 20 which contains a closed loop in the letter n.}
    \label{fig:nu}
\end{figure}

The high asymmetry score of this object is due to the fact that essential homology classes cannot be paired with other classes in a Wasserstein transport plan.  This is enforced by the definition of $d_{\Theta}((s,\textbf{Ord}),(t,\textbf{Rel})) = \infty$, which means it is always less costly to pair an essential interval to an ephemeral interval than to an interval from the ordinary or relative parts of the XPH diagram.  
This means the interval for the loop in the Modern No.~20 `n' cannot be matched with the corresponding non-essential class of the `u'. 

This example illustrates how our methodology is sensitive to alterations of objects that change their topology. While this is unsurprising for a topological method it does point to a limitation in our method to quantify approximate symmetries between objects when a small geometric change alters the essential homology classes without greatly altering a perceived symmetry, e.g., comparing a circle to a circle with a tiny segment removed. 

While these example research questions are somewhat quaint they demonstrate that the XPHT methodology can be used to quantify observed asymmetries, making it appropriate for researchers wishing to study symmetry in real systems such as the work in \cite{Biopaper} investigating which genes are responsible for the symmetry of leaves about the midvein.

\section{Conclusions}

We can see that the advantages of the method of symmetry quantification we have developed are that it requires little specification on the part of researchers, requiring only an image and a number of directions before outputting scores for all orthogonal transformations that leave the set of directions invariant. We can also gain good visual intuition for the asymmetry scores obtained by using polar plots to present these scores.

The limitations of this method are its sensitivity to changes in topology when comparing two objects, the current need for data to be in the form of a binary image, and how this restricts us to global shape symmetry under orthogonal transformations rather than localised partial symmetry. 

We should also mention that this method doesn't have objective units meaning it is more useful to compare two similar objects rather than interpreting the scores in isolation. For example, in Figure 1.a), the circle should be symmetric under all transformations and all asymmetry scores of the object are below 0.06. This means when the similar images in Figure 1.b-d) have transformations whose scores fall below 0.06, we can say with greater confidence that the object is symmetric under those transformations, at least when taken relative to Figure 1.a). Similarly, when the scores are above 0.06 in Figure 1.b-d) we can conclude that the object has only approximate symmetry under those transformations.

\subsection{Future Work}

The future directions for this work are to implement the XPHT for three-dimensional data so that we can analyse a broader class of objects. A notable example is protein structure where there are important relationships between symmetry and function and vast datasets of protein skeletons available for analysis. 
And although we can only analyse finite objects, the symmetry quantification method could be extended to analyse the presence of (approximate) translational isometry within an object. 

Another direction to explore is the analysis of a real-valued function by selecting various thresholds to create a set of sub-objects. These sub-objects could then be analysed using the XPHT and asymmetry scores.  It may also be possible to develop methods for partial symmetry detection within an object. 



\acknowledgments{
We would like to acknowledge the work of James Morgan in creating the R-package to calculate the XPHT of a binary image \cite{R}.}

\bibliographystyle{abbrv-doi}

\begin{thebibliography}{88} 

 \bibitem{networks} M. E. Aktas, E. Akbas, and  A. E. Fatmaoui.  Persistence homology of networks: methods and applications. \textit{Applied Network Science}, 4(61), 2019.

\bibitem{XPH} D. Cohen-Steiner, H. Edelsbrunner, and J. Harer. Extending Persistence Using Poincar\'e and Lefschetz Duality, \textit{Foundations of Computational Mathematics}, 9:79-103, 2009.

\bibitem{brain} L. Crawford, A. Monod, A. X. Chen, S. Mukherjee, and R. Rabad'an. Predicting Clinical Outcomes in Glioblastoma: An Application of Topological and Functional Data Analysis. \textit{Journal of the American Statistical Association}, 115:1139-1150, 2019.

\bibitem{DPC} V. De Silva, D. Morozov, and M. Vejdemo-Johansson. Dualities
in persistent (co) homology. \textit{Inverse Problems}, 27(12): 124003, 2011.

\bibitem{PHO}  H. Edelsbrunner, and J. Harer. Persistent homology - a survey. \textit{Discrete and Computational Geometry}, 453, 2008

\bibitem{TIA} P. Gandhi, M.-V. Ciocanel, K. Niklas, and A. T. Dawes. Identification of approximate symmetries in biological development. \textit{Philosophical Transactions. Series A, Mathematical, Physical and Engineering Sciences}, 379(2213):20200273, 2021.

\bibitem{SIG} N. J. Mitra, M. Pauly, M. Wand, and D. Ceylan. Symmetry in 3D Geometry: Extraction and Applications. \textit{Computer Graphics Forum}, 32(6), 2013.

\bibitem{R} J. Morgan. Source code for the Extended Persistent Homology Transform R Package. \textit{GitHub}, https://github.com/james-e-morgan/xpht, 2022.

\bibitem{PRST} J. Podolak, P. Shilane, A. Golovinskiy, S. Rusinkiewicz, and T. Funkhouser. A Planar-Reflective Symmetry Transform for 3D Shapes. \textit{ACM Transactions on Graphics (Proc. SIGGRAPH)}, 25(3), 2006.

\bibitem{bones} Y. Pritchard, A. Sharma, C. Clarkin, H. Ogden, S. Mahajan, and R. J. Sánchez-García. Persistent homology analysis distinguishes pathological bone microstructure in non-linear microscopy images. \textit{Scientific Reports}, 13(1):2522, 2023.

\bibitem{MPEG} T. Sikora. The MPEG-7 visual standard for content description-an overview. \textit{IEEE Transactions on circuits and systems for video technology}, 11(6):696–702, 2001.

\bibitem{leaf} P. F. B. Silva, A. R. S. Marcal, and R. M. Almeida da Silva. Evaluation of Features for Leaf Discrimination. \textit{International Conference on Image Analysis and Recognition}, 2013.

\bibitem{2011} D. M. Thomas, and V. Natarajan. Symmetry in Scalar Field Topology. \textit{IEEE Transactions on Visualization and Computer Graphics} 17(12):2035-2044, 2011

\bibitem{AEG} D. M. Thomas, and V. Natarajan. Detecting Symmetry in Scalar Fields Using Augmented Extremum Graphs. \textit{IEEE Transactions on Visualization and Computer Graphics}, 19(12):2663-2672, 2013.

\bibitem{PHT} K. Turner, S. Mukherjee, and D. M. Boyer. Persistent Homology Transform for Modelling Shapes and Surfaces. \textit{Information and Inference: A Journal of the IMA}, 3:310–344, 2014.

\bibitem{XPHT} K. Turner, V. Robins, and J. Morgan. The Extended Persistent Homology Transform of
Manifolds with Boundary. \textit{arXiv preprint arXiv:2208.14583}, 2022.

\bibitem{TI} G. V. Vstovsky. Transform information: A symmetry breaking measure. \textit{Foundations of Physics}, 27:1413-1444, 1997.

\bibitem{Biopaper} D. Wilson-Sánchez, S. Martínez-López, S. Navarro-Cartagena, S. Jover-Gil, and J. L. Micol. Members of the DEAL subfamily of the DUF1218 gene family are required for bilateral symmetry but not for dorsoventrality in Arabidopsis leaves. \textit{The New Phytologist}, 217(3):1307-1321, 2018.

\end{thebibliography}

\end{document}